\documentclass{amsart}
\usepackage{graphicx} % Required for inserting images
\usepackage[T1]{fontenc}
\usepackage[utf8]{inputenc}
\usepackage[english]{babel}
\usepackage{csquotes}
\usepackage{todonotes}
\usepackage{mathtools}

\usepackage[maxbibnames=99]{biblatex}
\renewbibmacro{in:}{%
        \ifentrytype{article}{}{\printtext{\bibstring{in}\intitlepunct}}}

\usepackage[hang,flushmargin]{footmisc}

\usepackage{amsmath}
\usepackage{amssymb}
\usepackage{amsthm}

\usepackage{abstract}
\addto\captionsenglish{}

\newtheorem{defn}{Definition}[section]
\newtheorem{theo}[defn]{Theorem}

\newtheorem{lemma}[defn]{Lemma}
\newtheorem{cor}[defn]{Corollary}

\newtheorem*{theov}{Theorem}

\newcommand{\N}{\mathbb{N}}

\newcommand{\Q}{\mathbb{Q}}
\newcommand{\R}{\mathbb{R}}

\newcommand{\func}[3]{#1:#2\rightarrow#3}

\newcommand{\Int}{\textup{int}}

\renewcommand{\epsilon}{\varepsilon}
\renewcommand{\theta}{\vartheta}

\title{Generic nonexpansive Hilbert space mappings}
\author[D.\ Ravasini]{Davide Ravasini $^1$}
\address[1]{Universit\"{a}t Leipzig, Mathematisches Institut, Augustusplatz 10, 04109 Leipzig, Germany}
\address{Charles University, Mathematical Institute, Sokolovsk\'{a} 49/83, 186 75 Prague, Czech Republic}
\email{davide.ravasini@matfyz.cuni.cz}
\author[D.\ K.\ Thimm]{Daylen K.\ Thimm $^2$}
\address[2]{Universit\"{a}t Innsbruck, Department of Mathematics, Technikerstra{\ss}e 13, 6020 Innsbruck, Austria}
\email{daylen.thimm@student.uibk.ac.at}

\begin{document}
\maketitle
\let\thefootnote\relax\footnote{\today \newline \indent\emph{2020 Mathematics Subject Classification:} 46C05, 54E52. \newline
}

\begin{abstract}
\noindent \textsc{Abstract}. We consider a closed convex set $C$ in a separable, infinite dimensional Hilbert space and endow the set $\mathcal{N}(C)$ of nonexpansive self-mappings on $C$ with the topology of pointwise convergence. We introduce the notion of a somewhat bounded set and establish a strong connection between this property and the existence of fixed points for the generic $f\in\mathcal{N}(C)$, in the sense of Baire categories. Namely, if $C$ is somewhat bounded, the generic nonexpansive mapping on $C$ admits a fixed point, whereas if $C$ is not somewhat bounded, the generic nonexpansive mapping on $C$ does not have any fixed points. This results in a topological 0--1 law: the set of all $f\in\mathcal{N}(C)$ with a fixed point is either meager or residual. We further prove that, generically, there are no fixed points in the interior of $C$ and, under additional geometric assumptions, we show the uniqueness of such fixed points for the generic $f\in\mathcal{N}(C)$ and the convergence of the iterates of $f$ to its fixed point.
\end{abstract}

\section{Introduction}
The question whether nonexpansive mappings have fixed points when defined on closed, convex subsets of Hilbert spaces and more general spaces has been investigated from various viewpoints. Let us recall that on a metric space $(X,\rho)$ a mapping $f\colon X\to X$ is nonexpansive if for all $x,y\in X$ the mapping satisfies
\[\rho\bigl(f(x),f(y)\bigr)\leq \rho(x,y).\]
A fundamental result is the classical theorem of Brouwer, which asserts that every continuous self-mapping of a bounded, closed, and convex subset of a Euclidean space has a fixed point. Notably, this theorem only works for the finite dimensional case. Generalizing it to infinite dimensions necessitates to only consider nonexpansive mappings and additional geometric properties, as demonstrated in the fixed point theorem of Browder--G\"{o}hde--Kirk~(see \cite{browder1965fixed, goebel1990topics, goebel1984uniform}). Finding a geometric restriction for a similar theorem to hold for all Lipschitz mappings is impossible: in every infinite dimensional Banach space there is a Lipschitz retraction from the ball to the sphere~\cite{benyamini1983spheres}, and the combination of this retraction with a suitable (also Lipschitz) reflection constitutes a Lipschitz mapping without a fixed point. Moreover, in some Banach spaces, e.g. $\mathcal{C}([0,1])$, there are nonexpansive mappings on closed, convex, and bounded subsets without a fixed point.\\

This brings up the following question:
\begin{quote}
	In infinite dimensions, how big is the set of nonexpansive mappings with a fixed point?
\end{quote}

\medskip
A pioneering result in this direction was found by de Blasi and Myjak~\cite{blasi1976convergence, blasi1989porosite}, who showed that on a bounded, closed, and convex subset of a Banach space, the generic nonexpansive mapping has a fixed point. Here generic is used to mean large in a topological sense. As the nonexpansive mappings with a fixed point form a large subset, one may wonder whether the strict contractions, i.e., Lipschitz mappings with Lipschitz constant $L<1$, are responsible for this, as they all possess a fixed point. It turns out that this is not the case, as de Blasi and Myjak had already shown in~\cite{blasi1976convergence} within the context of a Hilbert space that the set of strict contractions is negligible. Later it was shown that this result carries to Banach spaces in~\cite{bargetz2016sigma} and a large class of metric spaces in~\cite{bargetz2017porosity}. The type of mappings to be held responsible for the original result by de Blasi and Myjak have been found to be the Rakotch contractions, which all have fixed points. Recall that a Rakotch contraction is a mapping in which the distance of two images is not linearly bounded by the distance of the original points as in Lipschitz contractions, but by a factor that may depend on the distance. More precisely, a function $f$ on a metric space $(X, \rho)$ is a Rakotch contraction if there exists a decreasing function $\phi_f\colon [0,\infty) \to [0,1]$ with $\phi_f(t)<1$ for all $t>0$ such that for all $x,y \in X$ we have
\[\rho\bigl(f(x), f(y)\bigr)\leq \phi_f\bigl(\rho(x,y)\bigr)\rho(x,y).\]
Due to the Rakotch fixed point theorem~\cite{rakotch1962note}, all Rakotch contractions on a complete metric space have a fixed point. Furthermore, as shown by Reich and Zaslavski in~\cite{reich2000almost, reich2001set}, almost all nonexpansive mappings on a bounded, closed, and convex set are Rakotch contractions.\\

Whereas the set of nonexpansive functions on \emph{bounded} sets carries the topology of uniform convergence as a natural choice, in the \emph{unbounded} case, one has multiple topologies that seem natural to choose from, such as the topology of uniform convergence on bounded sets and the topology of pointwise convergence. Moreover, there are multiple natural metrics generating these topologies. In contrast to the bounded case, where there are many Rakotch contractions, Strobin was able to prove that with respect to the topology of uniform convergence on bounded sets the opposite holds: the set of Rakotch contractions is negligible~\cite{strobin2012some}. In fact, Strobin proved something more general. Recall that for $t\geq 0$ the modulus of continuity of $f$ is defined by
\[\omega_f(t)\coloneqq \sup\bigl\{\rho\bigl(f(x), f(y)\bigr)\colon x,y\in X, \rho(x,y)\leq t\bigr\}.\]
Given a concave and increasing function $\omega\colon [0,\infty)\to [0,\infty)$ with $\omega(0)=0$, Strobin considered the space $C_\omega(X)$ of all continuous self-mappings $f$ of $X$ satisfying $\omega_f(t)\leq\omega(t)$ for all $t\geq 0$ with respect to the topology of uniform convergence on bounded sets. He showed that if $X$ is an unbounded, closed, and convex subset of a Hilbert space, then the subset of all mappings satisfying $\omega_f(t)<\omega(t)$ for some $t>0$ is meager. For the case $\omega(t)=t$ this coincides with the previous statement about Rakotch contractions. Using a different metric, which generates a topology finer than the topology of uniform convergence on bounded sets, Reich and Zaslavski were able to prove in \cite{reich2016two} that the typical mapping is a Rakotch contraction on bounded subsets. As an important tool in Strobins proof, he uses the Gr\"{u}nbaum--Zarantonello extension theorem, which generalizes the Kirszbraun--Valentine theorem to uniformly continuous mappings.
Both these theorems are limited to Hilbert spaces. Overcoming this limitation, Bargetz, Reich, and the second author were recently able to generalize Strobin's result for nonexpansive mappings to hyperbolic spaces~\cite{brt_2023}. Shortly thereafter, the first author was able to extend this to arbitrary $\omega$~\cite{ravasini2024generic}.\\

However, in \cite{brt_2023} it is shown that with respect to the topology of uniform convergence on bounded sets the generic nonexpansive mapping maps a ball into itself and is a Rakotch contraction on \emph{bounded} sets. The Rakotch fixed point theorem therefore implies that the generic nonexpansive mapping has a unique fixed point.\\

Furthermore, in~\cite{brt_2023} Bargetz, Reich, and the second author looked at the problem in the unbounded setting from the perspective of the topology of pointwise convergence. In an effort to prove the generic presence of fixed points, they showed that for the generic nonexpansive mapping the strict inequality $\rho(f(x),f(y))<\rho(x,y)$ holds for residually many pairs of points $(x,y)\in X\times X$, thus hinting towards some kind of contractivity. However, we in fact disprove that the generic existence of fixed points is always guaranteed in the topology of pointwise convergence. This also provides an answer to the question asked at the end of \cite{ravasini2024generic}.\\

We focus on separable Hilbert spaces and, for the first time, we observe that depending upon what type of unbounded set is given, the generic nonexpansive mapping may or may not have a fixed point. Indeed, we are able to precisely classify on what sets the generic nonexpansive mapping does or does not have a fixed point. Since the sets where the generic nonexpansive mapping does have a fixed point exhibit the features of both bounded and possibly unbounded sets, we name them \emph{somewhat bounded} sets. We call all other sets totally unbounded. We are able to formulate this in two theorems that appear in the manuscript as Theorem~\ref{theo:sbfixed} and Theorem~\ref{theo:nofix}, respectively.% The result is a topological 0--1--Law:
%Change things here to include the two theorems and the corollary

\begin{theov} 
    Let $C$ be a closed, convex, somewhat bounded set in a separable, infinite dimensional Hilbert space $H$. Then there is a dense open set $\mathcal{R}\subseteq\mathcal{N}(C)$ such that every $f\in\mathcal{R}$ has a fixed point.
\end{theov}

\begin{theov}
    Let $C$ be a closed, convex set in a separable, infinite dimensional Hilbert space $H$. If $C$ is totally unbounded, then there is a residual set $\mathcal{R}\subseteq\mathcal{N}(C)$ such that every $f\in\mathcal{R}$ fails to have any fixed points.    
\end{theov}

\noindent Note that these two statements do not logically imply each other, as there exist sets which are neither meager nor residual. As our main result we directly obtain Theorem~\ref{theo:top01}, a topological 0--1 Law.

\begin{theov}[0--1 Law]
    Let $C$ be a closed, convex set in a separable, infinite-dimen\-sional Hilbert space $H$ and let $\mathcal{F}$ be the set of all $f\in\mathcal{N}(C)$ with a fixed point. Then either $\mathcal{F}$ is meager or it contains a dense open set.
\end{theov}

%\begin{quote}
	%Nonexpansive mappings on somewhat bounded sets generically have a fixed point, on all other sets they generically do not.
%\end{quote}
\medskip
As an additional result we show that the generic nonexpansive mapping on a somewhat bounded set has all fixed points on the boundary of the set. If the set is additionally locally uniformly rotund, we show that any sequence of iterates converges to the fixed point located on the boundary. Since our proofs rely on the Kirszbraun--Valentine extension theorem, the fixed point theorem of Browder--G\"{o}hde--Kirk and the classical machinery of Hilbert spaces, our results are limited to subsets of a Hilbert space.
\section{Setting and notation}
Given a metric space $(X,\rho)$, we denote by $\mathcal{N}(X)$ the space of nonexpansive mappings $\func{f}{X}{X}$. If $X$ is separable, we always assume $\mathcal{N}(X)$ to be endowed with the topology of pointwise convergence. This is the coarsest topology containing all open sets of the form $U_f(x,\epsilon)$, where $f\in\mathcal{N}(X)$, $x\in X$, $\epsilon>0$ and
\[ U_f(x,\epsilon)=\bigl\{g\in\mathcal{N}(X)\,:\,\rho\bigl(g(x),f(x)\bigr)<\epsilon\bigr\}. \]
If $X$ is complete and separable, the topology of pointwise convergence can be metrized via a complete metric as follows. We consider a dense, countable set $\Theta=\{\theta_n\}_{n=1}^\infty$ and define $\func{d_\Theta}{\mathcal{N}(X)\times\mathcal{N}(X)}{\R}$ by
\[ d_\Theta(f,g)=\sum_{n=1}^\infty 2^{-n}\frac{\rho\bigl(f(\theta_n),g(\theta_n)\bigr)}{1+\rho\bigl(f(\theta_n),g(\theta_n)\bigr)}. \]
The proof that $d_\Theta$ is indeed a complete metric on $\mathcal{N}(X)$ and that it induces the topology of pointwise convergence can be found in \cite{brt_2023}. We stress that the existence of a complete metric for $\mathcal{N}(X)$ makes the use of Baire categories meaningful, as we shall recall at the end of this section.

The following two lemmas will be needed later.
\begin{lemma}
    \label{lemma:1}
    Let $(X,\rho)$ be a metric space. For every $f\in\mathcal{N}(X)$, $x\in X$, $\epsilon>0$ and every positive integer $k$ the set
    \[ V_f(x,\epsilon,k)=\bigl\{g\in\mathcal{N}(X)\,:\,\rho\bigl(g^k(x),f^k(x)\bigr)<\epsilon\bigr\} \]
    is an open neighborhood of $f$ in $\mathcal{N}(X)$.
\end{lemma}
\begin{proof}
    We use induction on $k$. For $k=1$, note that $V_f(x,\epsilon,1)=U_f(x,\epsilon)$ for every $f\in\mathcal{N}(X)$, $x\in X$ and $\epsilon>0$. Let now $k$ be fixed and assume that $V_g(x,\delta,k)$ is open in $\mathcal{N}(X)$ for every $g\in\mathcal{N}(X)$, $x\in X$ and $\delta>0$. Given now $f\in\mathcal{N}(X)$, $x\in X$ and $\epsilon>0$, we have to show that $V_f(x,\epsilon,k+1)$ is open. Let $g\in V_f(x,\epsilon,k+1)$. The goal is to find an open neighborhood $W$ of $g$ such that $\rho(h^{k+1}(x),f^{k+1}(x))<\epsilon$ for every $h\in W$. Pick $\delta>0$ such that $2\delta+\rho(f^{k+1}(x),g^{k+1}(x))<\epsilon$ and set $W_1=U_g(g^k(x),\delta)$, $W_2=V_g(x,\delta,k)$. $W_1$ is open by the definition of the topology, whereas $W_2$ is open by the inductive assumption. Now set $W=W_1\cap W_2$. If $h\in W$ we have
    \begin{align*}
        \rho\bigl(h^{k+1}(x),f^{k+1}(x)\bigr) &\leq \rho\bigl(h(h^k(x)),h(g^k(x))\bigr)+\rho\bigl(h(g^k(x)),g(g^k(x))\bigr)+ \\
				      &\qquad +\rho\bigl(g^{k+1}(x),f^{k+1}(x)\bigr) < \\
				      &< \rho\bigl(h^k(x),g^k(x)\bigr)+\delta+\rho\bigl(g^{k+1}(x),f^{k+1}(x)\bigr)< \\
				      &<2\delta+\rho\bigl(f^{k+1}(x),g^{k+1}(x)\bigr)<\epsilon,
    \end{align*}
    which is what we wanted.
\end{proof}

\begin{lemma}
    \label{lemma:2}
    Let $(X,\rho)$ be a metric space, let $f$ be a nonexpansive mapping on $X$ and let $x\in X$ and $r>0$ be such that $\rho(f(x),x)>2r$. If $y\in X$ is a fixed point of $f$, then $\rho(x,y)>r$.
\end{lemma}
\begin{proof}
    We have
    \[ \rho(x,y)=\rho\bigl(x,f(y)\bigr)\geq\rho\bigl(x,f(x)\bigr)-\rho\bigl(f(x),f(y)\bigr)>2r-\rho(x,y), \]
    from which $\rho(x,y)>r$ follows.
\end{proof}

Throughout the paper we will deal with a single real, separable Hilbert space $H$ and its closed subspaces. We simply denote by $\langle x,y\rangle$ and $\|x\|$ the scalar product of two elements $x,y\in H$ and the norm of $x$ respectively. $B_H$ and $S_H$ stand for the closed unit ball and the unit sphere of $H$ respectively. Given a closed, convex set $C\subseteq H$, we denote by $\func{\pi_C}{H}{H}$ the nearest point projection onto $C$. In the case of a closed subspace $F\subseteq H$, we use the notation $\pi_F$ for the orthogonal projection onto $F$. This does not lead to any confusion, since the nearest point projection and the orthogonal projection coincide. Finally, given a subset $C\subseteq H$ and $x\in H$ we denote by $\textup{dist}(x,C)$ the distance of $x$ to $C$, i.e.,
\[ \textup{dist}(x,C)=\inf\{\|x-y\|\,:\,y\in C\}. \]

We recall that Hilbert spaces are uniformly convex Banach spaces. A Banach space $X$ is \emph{uniformly convex} if for every $\epsilon>0$ there is $\delta>0$ such that for every $x,y\in X$ with $\|x\|=\|y\|=1$ and $\|x-y\|>\delta$ we have $\|x+y\|\leq 2-\epsilon$. In uniformly convex Banach spaces, the following theorem holds (see \cite{browder1965fixed, goebel1984uniform, goebel1990topics}).
\begin{theo}[Browder--G\"{o}hde--Kirk]
    \label{thm:BGK}
    Let $C$ be a closed, convex, bounded set in a uniformly convex Banach space $X$. Then every nonexpansive mapping $\func{f}{C}{C}$ has a fixed point.
\end{theo}
\noindent Theorem \ref{thm:BGK} makes questions concerning the generic existence of fixed points rather trivial in the bounded case, which is the reason why our primary focus is on unbounded sets. Finally, we recall that in Hilbert spaces the Kirszbraun--Valentine Theorem on the extension of Lipschitz mappings is available (see, for instance, \cite{benlin}, Theorem 1.12 or \cite{reich2005fenchel}). This will be needed in the proof of Theorem \ref{theo:sbfixed}.
\begin{theo}[Kirszbraun--Valentine]
    \label{thm:kirszbraun}
    Let $H_1$ and $H_2$ be two Hilbert spaces, let $A\subseteq H_1$ be any subset and let $\func{f_0}{A}{H_2}$ be a Lipschitz mapping with Lipschitz constant $L$. Then there is a Lipschitz mapping $\func{f_1}{H_1}{H_2}$ with Lipschitz constant $L$ and such that $f_1(x)=f_0(x)$ for every $x\in A$.
\end{theo}

To conclude this preliminary section, we shall spend a few words on Baire category, since it plays a central role in the statement of all our results. A subset $N$ of a topological space $X$ is \emph{nowhere dense} if the closure of $N$ in $X$ has empty interior or, equivalently, if $X\setminus N$ contains a dense open set. A subset of $X$ is \emph{meager} if it is a countable union of nowhere dense sets. The complement of a meager set is often called \emph{comeager} or \emph{residual} and, by definition, it contains a countable intersection of dense open sets. The importance of meager and residual sets lies in the Baire category theorem, which asserts that a completely metrizable space is never meager in itself, i.e., it cannot be a countable union of some of its nowhere dense subsets. Since countable unions of meager sets are also meager, the Baire category theorem implies that a countable family of meager sets will never be able to cover the whole space. This conveys the intuitive idea that meager sets are small, whereas residual sets almost fill the entire space. The points of a residual set are usually called the \emph{generic} points of the space and a property which holds for every point in a residual set is usually said to be generic. For more on genericity in nonlinear analysis we refer the reader to \cite{reich2014genericity}.

\section{Somewhat bounded sets}
The central definition of our work is the following.
\begin{defn}
    A nonempty, non-singleton convex set $C$ in an infinite dimensional Hilbert space $H$ is \textbf{somewhat bounded} if there are $x_0\in C$, a finite dimensional subspace $F\subseteq H$ and $\alpha>0$ such that $\alpha B_F\subseteq C-x_0$ and $F^\perp\cap(C-x_0)$ is bounded. We say that $C$ is \textbf{totally unbounded} if it is not somewhat bounded.
\end{defn}
\noindent For example, the set $\{x\in H\,:\,\textup{dist}(x,F)\leq r\}$, where $F\subseteq H$ is a finite dimensional subspace and $r\geq 0$, is unbounded but also somewhat bounded.  $H$ itself is totally unbounded, as well as all its infinite dimensional subspaces.

Given a Hilbert space $H$, for every closed subspace $F\subseteq H$ and $\alpha,\beta>0$ we define the set
\[ D_F(\alpha,\beta)=\bigl\{ x\in H\,:\,\|\pi_{F^\perp}(x)\|-\beta\leq\alpha^{-1}\beta\|\pi_F(x)\| \bigr\}. \]
Note that, by definition, $(x+F^\perp)\cap D_F(\alpha,\beta)$ is bounded for every $x\in H$.
\begin{lemma}
\label{lemma:sb1}
Let $C$ be a convex set in an infinite dimensional Hilbert space $H$. Then the following are equivalent.
\begin{enumerate}
    \item $C$ is somewhat bounded.
    \item There are $z_0\in C$, a finite dimensional subspace $F\subseteq H$ and $\alpha,\beta>0$ such that $\alpha B_F\subseteq C-z_0\subseteq D_F(\alpha,\beta)$.
\end{enumerate}
\end{lemma}
\begin{proof}
(1)$\Rightarrow$(2). Since $C$ is somewhat bounded, we can find $z_0\in C$, a finite dimensional subspace $F\subseteq H$ and $\alpha>0$ such that $\alpha B_F\subseteq C-z_0$ and $F^\perp\cap(C-z_0)$ is bounded. Since $F^\perp\cap(C-z_0)$ is bounded, there must be $\beta>0$ such that $\beta S_{F^\perp}\cap(C-z_0)=\varnothing$. We want to show that $C-z_0\subseteq D_F(\alpha,\beta)$. Assume that this is not the case and pick $y\in(C-z_0)\setminus D_F(\alpha,\beta)$. Set
\[ \lambda=\frac{\beta}{\|\pi_{F^\perp}(y)\|-\beta} \]
and $w=-\lambda\pi_F(y)$. Note that
\[ \|w\|=\lambda\|\pi_F(y)\|=\frac{\beta\|\pi_F(y)\|}{\|\pi_{F^\perp}(y)\|-\beta}<\beta\beta^{-1}\alpha=\alpha, \]
thus $w\in\alpha B_F\subseteq C-z_0$. Moreover, if we set
\[ \mu=\frac{1}{1+\lambda}=\frac{\|\pi_{F^\perp}(y)\|-\beta}{\|\pi_{F^\perp}(y)\|}\in(0,1), \]
we obtain
\[ \mu w+(1-\mu)y=\frac{\beta\pi_{F^\perp}(y)}{\|\pi_{F^\perp}(y)\|}\in\beta S_{F^\perp}\cap(C-z_0). \]
This contradicts $\beta S_{F^\perp}\cap(C-z_0)=\varnothing$. Thus, we conclude that $C-z_0\subseteq D_F(\alpha,\beta)$, as claimed.

(2)$\Rightarrow$(1). This is obvious, since $F^\perp\cap(C-z_0)\subseteq F^\perp\cap D_F(\alpha,\beta)=\beta B_{F^\perp}$ is a bounded set.
\end{proof}
\begin{figure}
    \includegraphics[scale=0.8]{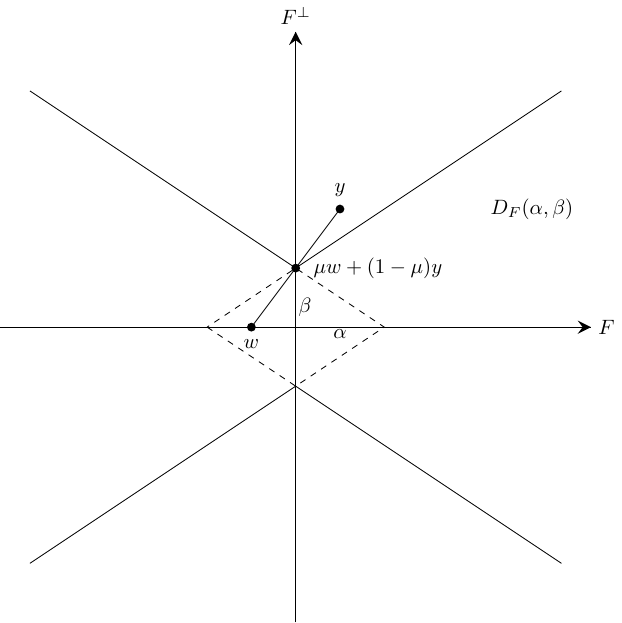}
    \caption{Visualization of the proof of Lemma~\ref{lemma:sb1}}
\end{figure}
\begin{lemma}
    \label{lemma:sb2}
    Assume that, for a closed, convex, unbounded set $C$ in a Hilbert space $H$, there are $\alpha,\beta>0$ and a subspace $F\subseteq H$ such that $\alpha B_F\subseteq C\subseteq D_F(\alpha,\beta)$. For every $\lambda,t>0$ define
    \[ E_F(\lambda,t)=\{x\in tS_F\,:\,\|\pi_C(x)\|\geq\lambda\}. \]
    Then, for every $\lambda>0$, there are $t\geq\lambda$ and $r,r'>0$ such that $r'>r$ and
    \[ r'B_H\cap C\subseteq E_F(\lambda,t)+rB_H. \]
\end{lemma}
\begin{proof}
    Let $\lambda>0$ be given. We first want to show that there is $t\geq\lambda$ such that $E_F(\lambda,t)$ is nonempty. To this end, select $t$ with
    \[ t-\frac{\beta}{\sqrt{\alpha^2+\beta^2}}(\alpha+t)\geq\lambda. \]
    Since $C$ is unbounded and $C\subseteq D_F(\alpha,\beta)$, there must $p\in C$ with $\|\pi_F(p)\|\geq t$. Define
    \[ v=\frac{t\pi_F(p)}{\|\pi_F(p)\|}\in tS_F,\quad w=-\frac{\alpha\pi_F(p)}{\|\pi_F(p)\|}\in\alpha S_F. \]
    The point on the line segment $\{w+\mu(p-w)\,:\,\mu\in[0,1]\}$ which minimizes the distance to $v$ is given by setting
    \[ \mu=\frac{\langle p-w,v-w\rangle}{{\|p-w\|}^2}. \]
    If we call $q$ such point, we have
    \[ {\|q-v\|}^2={(\alpha+t)}^2-\frac{{\langle p-w,v-w\rangle}^2}{{\|p-w\|}^2}. \]
    We now want to take a closer look at the right-hand side of the above equation. Using that $p\in D_F(\alpha,\beta)$ we get
    \[ \frac{\langle p-w,v-w\rangle}{\|p-w\|}=\|v-w\|\Bigl\langle\frac{p-w}{\|p-w\|},\frac{v-w}{\|v-w\|}\Bigr\rangle\geq(\alpha+t)\frac{\alpha}{\sqrt{\alpha^2+\beta^2}}. \]
    This implies
    \[ {\|q-v\|}^2\leq{(\alpha+t)}^2-\frac{\alpha^2}{\alpha^2+\beta^2}{(\alpha+t)}^2=\frac{\beta^2}{\alpha^2+\beta^2}{(\alpha+t)}^2. \]
    Since $q\in C$, we now have
    \[ \|v\|-\|\pi_C(v)\|\leq\|v-\pi_C(v)\|\leq\|v-q\|\leq\frac{\beta}{\sqrt{\alpha^2+\beta^2}}(\alpha+t), \]
    therefore
    \[ \|\pi_C(v)\|\geq\|v\|-\frac{\beta}{\sqrt{\alpha^2+\beta^2}}(\alpha+t)=t-\frac{\beta}{\sqrt{\alpha^2+\beta^2}}(\alpha+t)\geq\lambda, \]
    i.e.,\ $v\in E_F(\lambda,t)$.
    \begin{figure}
        \includegraphics[scale=0.8]{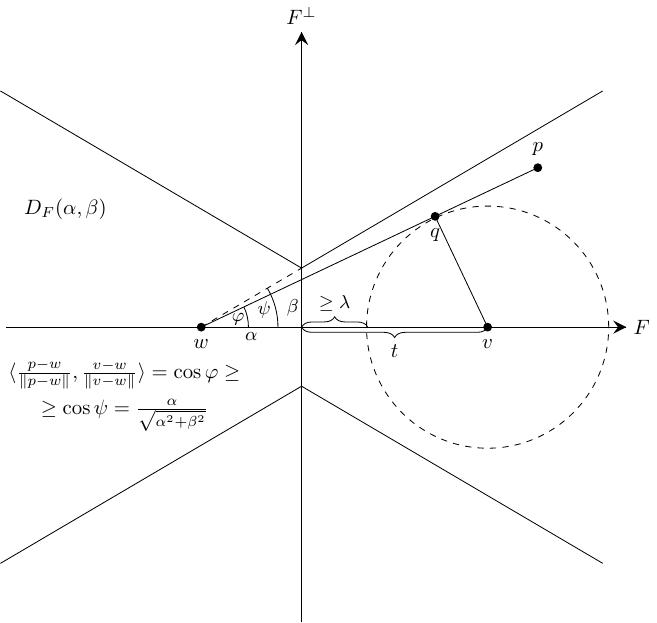}
        \caption{Visualization of the proof of Lemma~\ref{lemma:sb2} }
    \end{figure}
    We now set
    \[ r=\sqrt{{(2t)}^2+{(\beta+\alpha^{-1}\beta t)}^2} \]
    and $r'=\sqrt{r^2+t^2}>r$. Let $x\in C\cap r'B_H$. If $\|\pi_F(x)\|\leq t$, then $x\in \overline{B}(v,r)$, where $v\in E_F(\lambda,t)$ is the point we have determined before. We have indeed
    \[ {\|x-v\|}^2={\|\pi_F(x)-v\|}^2+{\|\pi_{F^\perp}(x)\|}^2\leq{(2t)}^2+{(\beta+\alpha^{-1}\beta t)}^2. \]
    If $\|\pi_F(x)\|>t$, then, as in the first part of the proof, we have
    \[ y=\frac{t\pi_F(x)}{\|\pi_F(x)\|}\in E_F(\lambda,t). \]
    Moreover,
    \begin{align*}
        {\|x-y\|}^2 &= {\|x\|}^2+{\|y\|}^2-2\langle x,y\rangle= \\
        &= {\|x\|}^2+t^2-2t\|\pi_F(x)\|\leq \\
	       &\leq r^2+2t^2-2t\|\pi_F(x)\|\leq \\
	       &\leq r^2+2t^2-2t^2= \\
	       &= r^2. \qedhere
    \end{align*}
\end{proof}

We are now ready to show the generic existence of fixed points for nonexpansive mappings on somewhat bounded sets.
\begin{theo}
    \label{theo:sbfixed}
    Let $C$ be a closed, convex, somewhat bounded set in a separable, infinite dimensional Hilbert space $H$. Then there is a dense open set $\mathcal{R}\subseteq\mathcal{N}(C)$ such that every $f\in\mathcal{R}$ has a fixed point.
\end{theo}
\begin{proof}
    First let us note that we may assume that $C$ is unbounded, as every nonexpansive mapping has a fixed point on bounded sets due to Theorem~\ref{thm:BGK}. By translating $C$ if needed, we may also assume that there are $\alpha,\beta>0$ and a finite dimensional subspace $F\subseteq H$ with $\alpha B_F\subseteq C\subseteq D_F(\alpha,\beta)$. Let $\Theta={\{\theta_n\}}_{n=1}^\infty$ be a dense set in $C$ and let $d_\Theta$ be the corresponding metric on $\mathcal{N}(C)$. To prove the existence of $\mathcal{R}$, let us fix some $f\in\mathcal{N}(C)$ and show that arbitrarily close to $f$ there is an open neighborhood of functions with a fixed point. Let us pick $\epsilon>0$ and find $n\in\N$ with $\sum_{j=n+1} ^\infty 2^{-j} < \epsilon$. Now set $M= \sup \{\|\theta_1\|, \dots,\|\theta_n\|, \|f(\theta_1)\|,\dots,\|f(\theta_n)\|\}$ and $\lambda=2M$. By Lemma \ref{lemma:sb2} we get $t\geq\lambda$ and $r', r >0$ with $r'>r$ such that
    \[r'B_H \cap C \subseteq E_F(\lambda, t) + rB_H, \]
    where $E_F(\lambda,t)$ is defined as in the statement of Lemma \ref{lemma:sb2}. Let us choose $\delta>0$ such that $r+\delta<r'$. Since $F$ is finite dimensional, $E_F(\lambda,t)$ is compact, therefore we can find a finite $\delta$-net $N\subseteq E_F(\lambda, t)$. Note that this implies
    \[ C\cap r'B_H \subseteq N + (r+\delta)B_H. \]
    Moreover, since the inequality $\|x-\pi_C(y)\|\leq\|x-y\|$ holds for every $x\in C$ and every $y\in H$, we also have $C\cap r'B_H\subseteq\pi_C(N)+(r+\delta)B_H$. Note further that $\|x\|\geq 2M$ for every $x\in\pi_C(N)$ by the definition of $E_F(\lambda,t)$. In particular, we have $\{\theta_1,\dots,\theta_n\}\cap\pi_C(N)=\varnothing$. Let us define the function
    \begin{align*}
        g_0: \{\theta_1,\dots,\theta_n\}\cup\pi_C(N) &\to C\\
						               x &\mapsto
						               \begin{cases}
							                 f(x) &\text{ if } x \in \{\theta_1,\dots,\theta_n\}\\
							                 0	 &\text{ if } x \in \pi_C(N)
						               \end{cases}
    \end{align*}
    and check that it is nonexpansive. Clearly, if $x,y \in \{\theta_1,\dots,\theta_n\}$ or $x,y\in \pi_C(N)$ we have that $\|g_0(x)-g_0(y)\|\leq \|x-y\|$, thus it only remains to show that the same holds for $x\in\{\theta_1,\dots,\theta_n\}$ and $y\in \pi_C(N)$. However, this is true as
    \[\|x-y\| \geq \|y\|-\|x\|\geq 2M - M = M\]
    and
    \[ \|g_0(x)-g_0(y)\| = \|f(x)\| \leq M \leq \|x-y\|. \]
    Now, Theorem \ref{thm:kirszbraun} yields a nonexpansive mapping $\func{g_1}{C}{H}$ which extends $g_0$, i.e., such that $g_1(x)=g_0(x)$ for every $x\in\{\theta_1,\dots,\theta_n\}\cup\pi_C(N)$. Finally, let us define $g=\pi_C\circ g_1\in\mathcal{N}(C)$ and note that $g$ also extends $g_0$. Furthermore, as $g$ and $f$ agree on $\{\theta_1,\dots,\theta_n\}$, we have that $d_\Theta(f,g)\leq \sum_{j=n+1}^\infty 2^{-j}<\epsilon$. We now want to exhibit an open neighborhood $U$ of $g$ in $\mathcal{N}(C)$ such that every $h\in U$ has a fixed point. Choose $\eta>0$ with $r+\delta+\eta<r'$. By the definition of the topology of pointwise convergence, the set
    \[ U=\{h\in\mathcal{N}(C)\,\colon\,\|h(x)\|<\eta\text{ for all }x\in\pi_C(N)\} \]
    is an open neighborhood of $g$ in $\mathcal{N}(C)$. For every $x\in r'B_H\cap C$ there is $y\in \pi_C(N)$ with $\|x-y\|\leq r+\delta$. If $h\in U$, we have that $\|h(y)\|<\eta$, hence
    \[ \|h(x)\|=\|h(x)-h(y)\|+\|h(y)\|< r+\delta+\eta < r'. \]
    The above line shows that $h(r'B_H\cap C) \subseteq r'B_H \cap C$ for every $h\in U$. Since every $h\in U$ maps a closed, convex, bounded subset of $C$ into itself, we can conclude by Theorem \ref{thm:BGK} that every $h\in U$ has a fixed point. Since $f$ and $\epsilon$ were chosen arbitrarily, this concludes the proof.
    %Hence, by Theorem~\ref{thm:BGK} we have that all $h\in\mathcal{R}'$ have a fixed point. As $f$ was chosen arbitrarily, we may construct such a set $\mathcal{R}'_f$ for every $f\in \mathcal{N}(C)$. Therefore, the set $\mathcal{R}\coloneqq \bigcup_{f\in\mathcal{N}(C)} \mathcal{R}'_f$ is a dense open set consisting only of functions with a fixed point.
\end{proof}

Although Theorem \ref{theo:sbfixed} guarantees the generic existence of fixed points in the case of a somewhat bounded set $C$, we shall now see that, generically, none of them lie in the interior of $C$. We start with a Lemma, which we also need in Section \ref{sec:iterates}.

\begin{lemma}
    \label{lemma:3}
    Let $C$ be a closed, convex, somewhat bounded set in a separable, infinite dimensional Hilbert space $H$. Then there is a residual set $\mathcal{R}\subseteq\mathcal{N}(C)$ such that for every $f\in\mathcal{R}$ and every $x\in C$ we have $\textup{dist}(f^k(x),\partial C)\to 0$ as $k\to\infty$.
\end{lemma}
\begin{proof}
    Let us assume without loss of generality that $0\in C$ and that $F\subseteq H$ is a finite dimensional subspace such that $C\cap(x+F^\perp)$ is bounded for every $x\in H$. Let $\Theta={\{\theta_n\}}_{n=1}^\infty$ be a dense set in $C$ and consider the corresponding metric $d_\Theta$. Given a positive integer $p$ and $r\in\Q\cap(0,1)$ define
    \[ R(p,r)=\bigl\{f\in\mathcal{N}(C)\,:\,\text{there is }k\text{ with }\textup{dist}\bigl(f^k(\theta_p),\partial C\bigr)<r\bigr\}. \]
    We now fix $p$ and $r$ arbitrarily and show that $R(p,r)$ contains a dense open set. To this end, we pick $f\in\mathcal{N}(C)$, $\delta>0$ and carry out the next part of the proof in two steps. \\

    \noindent\textbf{Step 1.} We want to find $h\in\mathcal{N}(C)$ with $d_\Theta(h,f)<\delta$ and a positive integer $k$ such that $h^k(\theta_p)\in\partial C$. Choose a positive integer $n$ with $\sum_{j=n+1}^\infty 2^{-j}<\delta/2$ and set $m=\max\{n,p\}$. Define
    \begin{align*}
        G &= F+\textup{span}\{\theta_1,\dots,\theta_m,f(\theta_1),\dots,f(\theta_m)\} \\
        D &= \textup{conv}\{\theta_1,\dots,\theta_m,f(\theta_1),\dots,f(\theta_m)\}
    \end{align*}
    Choose $y_0\in S_{G^\perp}$. Since $G^\perp\subseteq F^\perp$, we have that $C\cap(x+\R y_0)$ is compact for every $x\in D$, hence we can define the function $\func{\phi}{D}{[0,\infty)}$ by setting $\phi(x)=\max\{t\geq 0\,:\,x+ty_0\in C\}$ for every $x\in D$. The function $\phi$ is defined on a compact, convex subset of a finite dimensional Hilbert space and it is easily seen to be concave, thus it must be continuous, which means that it attains a maximum value $M\geq 0$. Let $\func{P}{H}{F}$ be the orthogonal projection onto $F$ and let $\func{Q}{H}{\R y_0}$ be the orthogonal projection onto $\R y_0$. Note that $\pi_D=\pi_D\circ P$. Define $g\in\mathcal{N}(C)$ as $g=\pi_D\circ f\circ\pi_D$ and $\func{h}{C}{C}$ by setting
    \[ h(x)=\pi_C\biggl(g(x)+Q(x)+\frac{\delta}{2} y_0\biggr) \]
    for every $x\in C$. We have to check that $h$ is nonexpansive. For every $x,y\in C$ we have
    \begin{align*}
         {\|h(x)-h(y)\|}^2 &\leq {\|g(x)-g(y)\|}^2+{\|Q(x)-Q(y)\|}^2= \\
 		  &= {\|(\pi_D\circ f\circ\pi_D\circ P)(x)-(\pi_D\circ f\circ\pi_D\circ P)(y)\|}^2+ \\
 		  &\quad +{\|Q(x)-Q(y)\|}^2\leq \\
		  &\leq {\|P(x)-P(y)\|}^2+{\|Q(x)-Q(y)\|}^2={\|x-y\|}^2,
    \end{align*}
    as we want. Next, we want $d_\Theta(h,f)<\delta$. Note that, for every $j\in\{1,\dots,m\}$, we have
    \begin{align*}
        \|h(\theta_j)-f(\theta_j)\| &= \biggl\|\pi_C\biggl(f(\theta_j)+\frac{\delta}{2}y_0\biggr)-\pi_C\bigl(f(\theta_j)\bigr)\biggr\|\leq \\
			    &\leq\biggl\|f(\theta_j)+\frac{\delta}{2}y_0-f(\theta_j)\biggr\|=\frac{\delta}{2},
    \end{align*}
    hence
    \[ d_\Theta(f,h)<\sum_{j=1}^m2^{-j}\|h(\theta_j)-f(\theta_j)\|+\frac{\delta}{2}<\frac{\delta}{2}+\frac{\delta}{2}=\delta. \]
    Finally, we want a positive integer $k$ with $h^k(\theta_p)\in\partial C$. Assume, by way of contradiction, that $h^k(\theta_p)\in\Int C$ for every positive integer $k$. Using that $\pi_C(x)\in\textup{int}C$ only if $\pi_C(x)=x$, it is not hard to show by induction on $k$ that the formula
    \[ h^k(\theta_p)=g^k(\theta_p)+\frac{k\delta}{2}y_0 \]
    holds for every $k\in\N$. Since $g^k(\theta_p)\in D$ for every $k$, we now have $k\leq 2\delta^{-1}M$ for every positive integer $k$, which is absurd. \\

    \noindent\textbf{Step 2.} We want to show that there is an open neighborhood $U\subseteq\mathcal{N}(C)$ of $h$ such that $\textup{dist}(u^k(\theta_p),\partial C)<r$ for every $u\in U$, where $k$ is the positive integer we have determined in Step 1. Using the notation of Lemma \ref{lemma:1}, this can be achieved by setting $U=V_h(\theta_p,r,k)$,
    which is indeed an open neighborhood of $h$ such that every $u\in U$ has the desired property. \\

    Step 1 and Step 2 prove that $R(p,q)$ contains a dense open set, as wished. Now define
    \[ \mathcal{R}=\bigcap_{r\in\Q\cap(0,1)}\bigcap_{p=1}^\infty R(p,r) \]
    and observe that $\mathcal{R}$ is residual. Let $f\in\mathcal{R}$, let $x\in C$ and let $\epsilon>0$. Find a positive integer $n$ with $\|x-\theta_n\|<\epsilon/2$. By our choice of $f$ we can find a positive integer $k$ with $\textup{dist}(f^k(\theta_n),\partial C)<\epsilon/2$. We now have
    \[ \textup{dist}\bigl(f^k(x),\partial C\bigr)\leq\|f^k(x)-f^k(\theta_n)\|+\textup{dist}\bigl(f^k(\theta_n),\partial C\bigr)\leq\|x-\theta_n\|+\frac{\epsilon}{2}<\epsilon . \]
    Since $x$ and $\epsilon$ are arbitrary, this shows that $\mathcal{R}$ is the set we are looking for.
\end{proof}

\begin{theo}
    \label{theo:sbbd}
    Let $C$ be a closed, convex, somewhat bounded set in a separable, infinite dimensional Hilbert space $H$. Then there is a residual set $\mathcal{R}\subseteq\mathcal{N}(C)$ such that, for every $f\in\mathcal{R}$, there are no fixed points of $f$ in $\Int C$.
\end{theo}
\begin{proof}
    Let $\mathcal{R}$ be the residual set given by Lemma \ref{lemma:3}. Let $f\in\mathcal{R}$ and assume that $x$ is a fixed point of $f$. Then
    \[ \textup{dist}(x,\partial C)=\textup{dist}\bigl(f^k(x),\partial C\bigr)\to 0 \]
    as $k\to\infty$. Thus, $\textup{dist}(x,\partial C)=0$. Since $\partial C$ is closed, this is possible only if $x\in\partial C$.
\end{proof}

\begin{cor}
    \label{cor:fixunique}
    Let $C$ be a closed, strictly convex, somewhat bounded set in a separable, infinite dimensional Hilbert space $H$. Then there exists a residual set $\mathcal{R}\subseteq\mathcal{N}(C)$ such that every $f\in\mathcal{R}$ has a unique fixed point.
\end{cor}
\begin{proof}
    Theorem \ref{theo:sbfixed} yields a dense open set $\mathcal{R}_1\subseteq\mathcal{N}(C)$ such that every $f\in\mathcal{R}_1$ has at least one fixed point, whereas Theorem \ref{theo:sbbd} yields a residual set $\mathcal{R}_2$ such that, for every $f\in\mathcal{R}_2$, all fixed points of $f$ lie in $\partial C$. Set $\mathcal{R}=\mathcal{R}_1\cap\mathcal{R}_2$ and note that $\mathcal{R}$ is residual. Note that, for every $f\in\mathcal{N}(C)$, the set of fixed points of $f$ is convex because the norm of $H$ is strictly convex. If $f\in\mathcal{R}$, then the set of fixed points of $f$ is nonempty and contained in $\partial C$. Since $\partial C$ does not admit any nontrivial convex subsets, we deduce that the fixed point of $f$ must be unique.
\end{proof}

\section{Totally unbounded sets}
In the present section we aim to show that the situation is completely opposite on totally unbounded sets, in that the generic nonexpansive mapping does not admit any fixed points. We can already formulate and prove the main theorem. The proof adapts and makes use of the same iteration technique we have already seen in Lemma \ref{lemma:3}.
\begin{theo}
\label{theo:nofix}
Let $C$ be a closed, convex set in a separable, infinite dimensional Hilbert space $H$. If $C$ is totally unbounded, then there is a residual set $\mathcal{R}\subseteq\mathcal{N}(C)$ such that every $f\in\mathcal{R}$ fails to have any fixed points.
\end{theo}
\begin{proof}
Without loss of generality, we may assume $0\in C$. Let $\Theta={\{\theta_n\}}_{n=1}^\infty$ be a dense subset in $C$ and define the corresponding metric $d_\Theta$. For every positive integer $r$, define $R(r)$ as the set of all $f\in\mathcal{N}(C)$ such that no fixed point of $f$ lies in $rB_H$. We want to show that $R(r)$ contains a dense open set for every $r$. To this end, we consider $f\in\mathcal{N}(C)$ and $\epsilon>0$ and divide the next part of the proof in two steps. \\

\noindent\textbf{Step 1.} We want to find $h\in\mathcal{N}(C)$ and a positive integer $k$ such that $d_\Theta(f,h)<\epsilon$ and $\|h^k(0)\|\geq 3r$. Choose a positive integer $n$ with $\sum_{j=n+1}^\infty 2^{-j}<\epsilon/2$ and set
\[ F = \textup{span}\{\theta_1,\dots,\theta_n,f(\theta_1),\dots,f(\theta_n)\}. \]
Note that $C\cap F$ has nonempty relative interior in $F$, which implies that there are $x_0\in C\cap F$ and $\alpha>0$ with $\alpha B_F\subseteq C-x_0$. Set
\[ \rho=\max\{\|\theta_1-x_0\|,\dots,\|\theta_n-x_0\|,\|f(\theta_1)-x_0\|,\dots,\|f(\theta_n)-x_0\|\} \]
and pick any positive $\delta<\min\{1,{(8\rho)}^{-1}\epsilon\}$. The assumption that $C$ is totally unbounded implies that $C\cap(x_0+F^\perp)$ is unbounded, hence we can find $y_0\in S_{F^\perp}$ such that $x_0+3r\delta^{-1} y_0\in C$. Define
\[ D = \delta x_0+(1-\delta)\cdot\textup{conv}\{\theta_1,\dots,\theta_n,f(\theta_1),\dots,f(\theta_n)\}. \]
It follows from the convexity of $C$ that \[ \delta x_0+(1-\delta)C+[0,3r]y_0\subseteq C. \]
Indeed, for every $x\in C$ and $t\in[0,3r]$ one has
\begin{align*}
\delta x_0+(1-\delta)x+ty_0 &\in (1-\delta)C+[\delta x_0,\delta x_0+3r y_0]= \\
                            &= (1-\delta)C+\delta[x_0,x_0+3r\delta^{-1}y_0]\subseteq(1-\delta)C+\delta C=C.
\end{align*}
In particular, one has $D+[0,3r]y_0\subseteq C$. Let $\func{P}{H}{F}$ and $\func{Q}{H}{\R y_0}$ be the orthogonal projections onto $F$ and $\R y_0$ respectively. Pick a positive integer $k>12r\epsilon^{-1}$ and define $g\in\mathcal{N}(C)$ as $g=\pi_D\circ f\circ\pi_D$ and $\func{h}{C}{C}$ by setting
\[ h(x)=\pi_C\bigl(g(x)+Q(x)+3rk^{-1}y_0\bigr) \]
for every $x\in C$. One checks that $h$ is nonexpansive similarly to the proof of Lemma \ref{lemma:3}. To see that $d_\Theta(f,h)<\epsilon$, note that, for every $j\in\{1,\dots,n\}$, one has
\begin{align*}
\|h(\theta_j)-f(\theta_j)\| &= \|\pi_D\circ f\circ\pi_D(\theta_j)+3rk^{-1}y_0-f(\theta_j)\|\leq \\
			    &\leq \|\pi_D\circ f\circ\pi_D(\theta_j)-\pi_D\circ f(\theta_j)\|+ \\
			    &\qquad+\|\pi_D\circ f(\theta_j)-f(\theta_j)\|+3r k^{-1}\leq \\
			    &\leq\|\pi_D(\theta_j)-\theta_j\|+\|\pi_D\circ f(\theta_j)-f(\theta_j)\|+3rk^{-1}\leq \\
			    &\leq\|\delta x_0+(1-\delta)\theta_j-\theta_j\|+ \\
              &\qquad +\|\delta x_0+(1-\delta)f(\theta_j)-f(\theta_j)\|+3rk^{-1}< \\
			    &< 2\delta\rho+\frac{\epsilon}{4}.
\end{align*}
Thus, by our choice of $\delta$ we obtain $\|h(\theta_j)-f(\theta_j)\|<\epsilon/2$ for every $j\in\{1,\dots,n\}$. It follows that
\begin{align*}
d_\Theta(h,f) &\leq \sum_{j=1}^n2^{-j}\|h(\theta_j)-f(\theta_j)\|+\sum_{j=n+1}^\infty 2^{-j}< \\
	      &< \frac{\epsilon}{2}\sum_{j=1}^\infty 2^{-j}+\sum_{j=n+1}^\infty 2^{-j}\leq\frac{\epsilon}{2}+\frac{\epsilon}{2}=\epsilon.
\end{align*}
Finally, using the inclusion $D+[0,3r]y_0\subseteq C$, it is not hard to check by induction that the formula $h^j(x)=g^j(x)+3rk^{-1}jy_0$ holds for every $x\in C\cap F$ and every $j\in\{1,\dots,k\}$. In particular, we have $h^k(0)=g^k(0)+3ry_0$. Since $g^k(0)$ and $y_0$ are orthogonal to each other, this implies
\[ {\|h^k(0)\|}^2={\|g^k(0)\|}^2+9r^2\geq 9r^2, \]
which is what we wanted. \\

\noindent\textbf{Step 2.}  We want to find a neighborhood $U$ of $h$ in $\mathcal{N}(C)$ such that every $u\in U$ does not have any fixed points in $rB_H$. Using the notation of Lemma \ref{lemma:1}, set $U=V_h(0,r,k)$, where $k$ is the positive integer we found in Step 1. If $x\in C$ is a fixed point of $u\in U$, then it is also a fixed point of $u^k$. Note that
\[ \|u^k(0)\|\geq\|h^k(0)\|-\|u^k(0)-h^k(0)\|>3r-r=2r, \]
therefore it follows from Lemma \ref{lemma:2} that $\|x\|>r$. \\

Step 1 and Step 2 show that $R(r)$ contains a dense open set for every $r$, as desired. Now set
\[ \mathcal{R}=\bigcap_{r=1}^\infty R(r), \]
and note that $\mathcal{R}$ is residual. If $f\in\mathcal{R}$ has a fixed point $x$, then $\|x\|>r$ for every positive integer $r$, which is impossible. Hence, $f$ does not have any fixed points.
\end{proof}

If we combine Theorem \ref{theo:sbfixed} and Theorem \ref{theo:nofix}, we obtain the topological 0--1 law we mentioned in the introduction.
\begin{theo}[0--1 Law]\label{theo:top01}
    Let $C$ be a closed, convex set in a separable, infinite dimensional Hilbert space $H$ and let $\mathcal{F}$ be the set of all $f\in\mathcal{N}(C)$ with a fixed point. Then either $\mathcal{F}$ is meager or it contains a dense open set.
\end{theo}

\section{Convergence of the iterates}
\label{sec:iterates}
We dedicate this last section to the study of the iterates $f^k$, $k\in\N$, of the generic nonexpansive mapping $f\in\mathcal{N}(C)$, where $C$ is a closed and convex set in a separable Hilbert space. Our convergence result will be obtained under a further geometric assumption on $C$, which is usually called \emph{local uniform rotundity}. We say that $C$ is \emph{locally uniformly rotund}, or LUR, if $\partial C\ne\varnothing$ and for every $x\in C$ and $\epsilon>0$ there is $\delta>0$ such that
\[ \textup{dist}\biggl(\frac{x+y}{2},\partial C\biggr)>\epsilon \]
whenever $y\in C$ and $\|x-y\|>\delta$. Note that LUR sets are always strictly convex and have therefore nonempty interior. Alternatively, one can define LUR sets via their \emph{modulus of convexity}. For every $x\in C$, set $\Delta_C(x)=\sup\{\|y-x\|\,:\,y\in C\}$. Now, for every $x\in C$ and every $\epsilon\in(0,\Delta_C(x))$, define
\[ \delta_C(x,\epsilon)=\inf\biggl\{\textup{dist}\biggl(\frac{x+y}{2},\partial C\biggr)\,:\,y\in C,\,\|y-x\|\geq\epsilon\biggr\}. \]
$C$ is locally uniformly rotund precisely when $\delta_C(x,\epsilon)>0$ for every $x\in C$ and every $\epsilon\in(0,\Delta_C(x))$. LUR sets are a natural generalization of LUR norms and have been considered, for example, in \cite{cvz_1996} and in \cite{berves_2023}, where the authors investigate this property in connection with existence of extensions for quasiconvex functions. We will use the fact that the modulus of convexity $\delta_C(x,\epsilon)$ is a continuous function of $\epsilon$ once $x$ is fixed.

\begin{lemma}
    \label{lemma:iterates1}
    Let $C$ be an LUR set in a Hilbert space $H$ and let $x,y,z$ be points in $C$ with $y\in\partial C$. If $r,s>0$ are such that $r\leq\|x-y\|\leq s$, then
    \[ \Bigl\langle\frac{z-y}{\|z-y\|},\frac{x-y}{\|x-y\|}\Bigr\rangle\geq-\frac{s}{\sqrt{s^2+4{\delta_C(x,r)}^2}}. \]
\end{lemma}
\begin{proof}
    We may assume without loss of generality that $y=0$. Define $F=\textup{span}\{x,z\}$ and note that, since $C$ is strictly convex, $x$ and $z$ cannot be linearly dependent, hence $\dim F=2$. Let $w\in S_F$ be such that $\langle x,w\rangle=0$ and, for every $t\in\R$, define
    \[ p_t=\frac{1}{2}x+tw. \]
    Set $\gamma=\delta_C(x,r)$. Observe that, since $\textup{dist}(x/2,\partial C)\geq\gamma$, we have that $p_t\in\Int C$ for every $t\in(-\gamma,\gamma)$. Thus, for every $t\in(-\gamma,\gamma)$, the half line $\{sp_t\,:\,s\leq 0\}$ does not contain any points of $C$ besides $0$. In particular, $C\cap F$ is contained in the cone
    \[ Q=\{\lambda x+\mu w\,:\,\lambda,\mu\in\R,\,2\gamma\lambda\geq-|\mu|\}. \]
    Let $\lambda, \mu$ be such that $z=\lambda x+\mu w$. If $\lambda\geq 0$ the statement is obvious. If $\lambda <0$, we have ${(\mu/\lambda)}^2\geq 4\gamma^2$. This implies
    \begin{align*}
        \Bigl\langle\frac{z}{\|z\|},\frac{x}{\|x\|}\Bigr\rangle &= \frac{\lambda\|x\|}{\sqrt{\lambda^2{\|x\|}^2+\mu^2}}=-\frac{\|x\|}{\sqrt{{\|x\|}^2+{(\mu/\lambda)}^2}}\geq \\
        & \geq-\frac{\|x\|}{\sqrt{{\|x\|}^2+4\gamma^2}}\geq-\frac{s}{\sqrt{s^2+4\gamma^2}},
    \end{align*}
    as claimed.
\end{proof}

\noindent\emph{Remark:} On a geometric level, Lemma \ref{lemma:iterates1} asserts that the angle $\varphi$ between the vectors $x-y$ and $z-y$ must be bounded away from $\pi$ by a quantity which depends on $r$ and $s$ and on the modulus of convexity of $C$ at $x$. More precisely, we have
\[ \varphi\leq\arccos\left(-\frac{s}{\sqrt{s^2+4{\delta_C(x,r)}^2}}\right)=\pi-\arccos\left(\frac{s}{\sqrt{s^2+4{\delta_C(x,r)}^2}}\right). \]

\begin{lemma}
    \label{lemma:iterates2}
    Let $C$ be an LUR set in a Hilbert space $H$ and let $\func{f}{C}{C}$ be a nonexpansive mapping such that $\textup{dist}(f^k(x),\partial C)\to 0$ for every $x\in C$ as $k\to\infty$. Moreover, let $x_0\in C$ be a fixed point of $f$. Then there is a continuous function $\func{\alpha}{(0,\Delta_C(x_0))}{(0,1)}$ such that for every $r\in(0,\Delta_C(x_0))$ and for every $x\in C$ with $\|x-x_0\|=r$ there is a positive integer $k$ with $\|f^k(x)-x_0\|\leq\alpha(r)r$.
\end{lemma}
\begin{proof}
    Let $\func{\delta}{[0,\Delta_C(x_0))}{[0,\infty)}$ be the modulus of convexity of $C$ at the point $p$. We introduce the function $\func{\beta}{(0,\Delta_C(x_0))}{(0,1)}$ given by
    \[ \beta(r)=\frac{5r}{\sqrt{25r^2+256{\delta(r/8)}^2}} \]
    and observe that
    \[ \sqrt{\frac{1+\beta(r)}{2}}<1 \]
    holds for every $r\in(0,\Delta_C(x_0))$. Therefore, we can find a continuous function $\func{\epsilon}{(0,\Delta_C(x_0))}{(0,1/4)}$ with
    \[ \sqrt{\frac{1+\beta(r)}{2}}\cdot\bigl(1+\epsilon(r)\bigr)<1. \]
    We claim that the function $\alpha$ we are looking for can be defined by
    \[ \alpha(r)=\max\Biggl\{\frac{3}{4},\,\sqrt{\frac{1+\beta(r)}{2}}\cdot\bigl(1+\epsilon(r)\bigr)\Biggr\}. \]
    Fix $r\in(0,\Delta_C(x_0))$, let $x\in C$ with $\|x-x_0\|=r$ and let $y$ be the midpoint of the segment $[x_0,x]$. By our assumption on $f$ there is a positive integer $k$ such that $\textup{dist}(f^k(y),\partial C)<r\epsilon(r)/4$, which allows us to choose a point $z\in\partial C$ with $\|f^k(y)-z\|<r\epsilon(r)/2$. Since $f$ is nonexpansive we now get
    \begin{align*}
        \|x_0-z\| \leq \|x_0-f^k(y)\|+\|f^k(y)-z\| &\leq\bigl(1+\epsilon(r)\bigr)\frac{r}{2}, \\
        \|f^k(x)-z\| \leq \|f^k(x)-f^k(y)\|+\|f^k(y)-z\| &\leq\bigl(1+\epsilon(r)\bigr)\frac{r}{2}.
    \end{align*}
    In particular, we have $\|x_0-z\|\leq 5r/8$. If $\|f^k(y)-x_0\|<r/4$, we get
    \[ \|f^k(x)-x_0\|\leq\|f^k(x)-f^k(y)\|+\|f^k(y)-x_0\|<\frac{3}{4}r\leq\alpha(r)r \]
    and the lemma holds. If $\|f^k(y)-x_0\|\geq r/4$, we have
    \[ \|x_0-z\|\geq\|f^k(y)-x_0\|-\|f^k(y)-z\|\geq\frac{r}{4}-\frac{r\epsilon(r)}{2}=\frac{r}{4}\bigl(1-2\epsilon(r)\bigr)\geq\frac{r}{8}. \]
    Now we want to apply Lemma \ref{lemma:iterates1} to the points $x_0$, $z$ and $f^k(x)$ and the bound $r/8\leq\|x_0-z\|\leq 5r/8$, thus obtaining
    \[ \Bigl\langle\frac{f^k(x)-z}{\|f^k(x)-z\|},\frac{x_0-z}{\|x_0-z\|}\Bigr\rangle\geq-\frac{5r/8}{\sqrt{{{(5r/8)}^2+4{\delta(r/8)}^2}}}=-\beta(r). \]
    Therefore,
    \begin{align*}
        {\|f^k(x)-x_0\|}^2 &= {\|f^k(x)-z\|}^2+{\|x_0-z\|}^2 \\ &\qquad-2\|f^k(x)-z\|\cdot\|x_0-z\|\cdot\Bigl\langle\frac{f^k(x)-z}{\|f^k(x)-z\|},\frac{x_0-z}{\|x_0-z\|}\Bigr\rangle\leq \\
        &\leq {\bigl(1+\epsilon(r)\bigr)}^2\cdot\frac{r^2}{4}+{\bigl(1+\epsilon(r)\bigr)}^2\cdot\frac{r^2}{4} \\
        &\qquad+2{\bigl(1+\epsilon(r)\bigr)}^2\cdot\frac{r^2}{4}\cdot\beta(r)= \\
        &= \frac{1}{2}\bigl(1+\beta(r)\bigr)\cdot{\bigl(1+\epsilon(r)\bigr)}^2r^2\leq{\alpha(r)}^2r^2,
    \end{align*}
    which again confirms our claim.
\end{proof}

\begin{theo}
    Let $C$ be a closed, somewhat bounded LUR set in a separable, infinite dimensional Hilbert space. Then there exists a residual set $\mathcal{R}\subseteq\mathcal{N}(C)$ such that every $f\in\mathcal{R}$ has a unique fixed point and the sequence of iterates ${(f^k(x))}_{k\in\N}$ converges to the fixed point of $f$ for every $x\in C$.
\end{theo}
\begin{proof}
    Theorem \ref{theo:sbfixed} yields a dense open set $\mathcal{R}_1\subseteq\mathcal{N}(C)$ such that every $f\in\mathcal{R}_1$ has a fixed point, whereas Lemma \ref{lemma:3} gives us a residual set $\mathcal{R}_2\subseteq\mathcal{N}(C)$ such that for every $f\in\mathcal{R}_2$ and every $x\in C$ we have $\textup{dist}(f^k(x),\partial C)\to 0$ as $k\to\infty$. Set $\mathcal{R}=\mathcal{R}_1\cap\mathcal{R}_2$. We have already concluded in Theorem \ref{theo:sbbd} and the subsequent Corollary \ref{cor:fixunique} that every $f\in\mathcal{R}$ has a unique fixed point, which in addition lies on $\partial C$. Let $f\in\mathcal{R}$ and let $x_0$ be its unique fixed point. Since $f$ satisfies the assumptions of Lemma \ref{lemma:iterates2}, we get a continuous function $\func{\alpha}{(0,\Delta_C(x_0))}{[0,1)}$ such that for every $x\in C$ there is a positive integer $k$ with $\|f^k(x)-x_0\|\leq\alpha(\|x-x_0\|)\|x-x_0\|$. Let $x\in C$, set $r_k=\|f^k(x)-x_0\|$ for every $k\in\N$ and note that the sequence ${(r_k)}_{k\in\N}$ is nonincreasing and has therefore a limit and infimum $r$. Assume that $r\ne 0$. Since $\alpha$ is continuous, we can find $\delta>0$ such that $\alpha(t)t<r$ for every $t\in[r,r+\delta)$. Let $k$ be such that $r_k<r+\delta$. Lemma \ref{lemma:iterates2} now yields $k'>k$ such that $r_{k'}\leq\alpha(r_k)r_k<r$. This contradicts the fact that $r$ is the infimum, hence $r=0$.
\end{proof}

\noindent\emph{Remark:} It was pointed out to the authors that the previous theorem can also be proven by combining our Theorem \ref{theo:sbbd} with Theorem 5.1 from \cite{Grzesik2019}, as in this theorem the boundedness of the set is only used to conclude the existence of a fixed point, which is already available in our case.

\section{Open questions}
For a closed, convex, somewhat bounded set $C$ in a separable Hilbert space $H$, the following questions are still open.
\begin{enumerate}
    \item Does the generic mapping $f\in\mathcal{N}(C)$ have a unique fixed point even if $C$ is not strictly convex?
    \item Do the iterates $f^k(x)$ of the generic mapping $f\in\mathcal{N}(C)$ converge to a fixed point of $f$ for every $x\in C$ even if $C$ is not an LUR set?
\end{enumerate}

However, the natural research direction concerns the possibility to extend our results to the Banach space setting. In what follows, $C$ is a closed, convex set in a separable, infinite-dimensional Banach space.
\begin{enumerate}
    \setcounter{enumi}{2}
    \item Does it hold that the generic mapping $f\in\mathcal{N}(C)$ fails to have any fixed points in $\Int C$?
    \item Let $\mathcal{F}$ be the set of all $f\in\mathcal{N}(C)$ with a fixed point. Is $\mathcal{F}$ either meager or residual?
    \item Let $\mathcal{F}$ be the set of all $f\in\mathcal{N}(C)$ which map a closed, convex, bounded subset of $C$ to itself. Is $\mathcal{F}$ either meager or residual?
    \item What can be said about the uniqueness of fixed points and the convergence of the iterates for the generic $f\in\mathcal{N}(C)$ in case the set of all $f\in\mathcal{N}(C)$ with a fixed point is residual?
\end{enumerate}

\section*{Acknowledgments}
 We thank Prof.\ Simeon Reich for his valuable suggestions and the referee for their careful reading of the manuscript. The second author would like to thank his supervisor Eva Kopeck\'{a} for her encouragement regarding this research. The second author was supported by the Doctoral Scolarship of the University of Innsbruck. He also received funding from the State of Tyrol (TNF).

\printbibliography

\end{document}